\documentclass[letterpaper, 10 pt, conference]{ieeeconf}  

\IEEEoverridecommandlockouts                              

\usepackage{amsmath,mathtools}
\usepackage{amssymb}
\usepackage{times}
\usepackage{multirow}
\usepackage{graphicx}
\usepackage[dvipsnames,table,xcdraw]{xcolor}
\usepackage{tikz}
\usepackage{subcaption}
\usepackage[normalem]{ulem}
\usepackage{algorithmic,algorithm}
\usepackage[noadjust]{cite}
\usetikzlibrary{arrows}
\usetikzlibrary{calc}

\newcommand{\Xs}{\mathcal{X}}

\newcommand{\Sem}{\textbf{\emph{s}}}
\newcommand{\Aem}{\textbf{\emph{a}}}
\newcommand{\Iem}{\textbf{\emph{i}}}
\newcommand{\Rem}{\textbf{\emph{r}}}
\newcommand{\Dem}{\textbf{\emph{d}}}

\newcommand*\diff{\mathop{}\!\mathrm{d}}

\newcommand{\yonote}[1]{\textcolor{red}{[YO: #1]}}

\def\Re{\mathbb{R}}

\def\argmin{\mathop{\text{\rm arg\,min}}}

\def\Sec#1{Sec.~\ref{#1}}

\def\notes#1{\marginpar{\tiny #1}\typeout{Notes!
Notes!
Notes!
}}
\renewcommand{\notes}[1]{\typeout{notes!}}

\def\Re{\field{R}}

\def\Sec#1{Sec.~\ref{#1}}

\def\Sec#1{Sec~\ref{#1}}

\def\E{{\sf E}}

\def\Sec#1{Sec.~\ref{#1}}

\newtheorem{assumption}{Assumption}
\newtheorem{definition}{Definition}

\newtheorem{remark}{Remark}
\newtheorem{proposition}{Proposition}

\def\beq{\begin{eqnarray}} 
\def\bc{\begin{center}} 
\def\be{\begin{enumerate}}
\def\bi{\begin{itemize}} 
\def\bs{\begin{small}}
\def\bS{\begin{slide}}
\def\ec{\end{center}} 
\def\ee{\end{enumerate}}
\def\ei{\end{itemize}}
\def\es{\end{small}}
\def\eS{\end{slide}}
\def\eeq{\end{eqnarray}}

\newcommand{\newP}[1]{\medskip\noindent{\bf #1:}}

\newcommand{\ud}{\,\mathrm{d}}

\def\Re{\mathbb{R}}
\def\E{{\sf E}}

\def\argmin{\mathop{\text{\rm arg\,min}}}

\def\Sec#1{Sec.~\ref{#1}}

\renewcommand{\Re}{\mathbb{R}}

\def\ls{\lambda^{\text{\tiny SA}}}
\def\la{\lambda^{\text{\tiny AI}}}
\def\lar{\lambda^{\text{\tiny AR}}}
\def\lr{\lambda^{\text{\tiny IR}}}
\def\ld{\lambda^{\text{\tiny ID}}}
\def\be{\lambda^{\text{\tiny SA}}\b_t}

\def\beb{\lambda^{\text{\tiny SA}}\bar{\b}}

\def\at{a^{\text{\tiny thresh}}}

\title{\LARGE \bf
How does a Rational Agent Act in an Epidemic?}

\author{S. Yagiz Olmez, Shubham Aggarwal, Jin Won Kim, Erik Miehling
  \\[5pt] Tamer Ba{\c s}ar, Matthew West, and Prashant G. Mehta\\[5pt]
\thanks{Research supported in part by the C3.ai Digital Transformation Institute sponsored by C3.ai Inc. and the Microsoft Corporation and in part by the National Science Foundation grants NSF-ECCS 20-32321 and NSF-CMMI 1761622.}
\thanks{	S. Y. Olmez, S.Aggarwal, J. W. Kim, and P. G. Mehta are with the Coordinated Science Laboratory and the Department of Mechanical Science and Engineering at the University of Illinois at Urbana-Champaign (UIUC); E. Miehling and T. Ba{\c s}ar are with the Coordinated Science Laboratory and the Department of Electrical and Computer Engineering at UIUC; M. West is with the Department of Mechanical Science and Engineering at UIUC.}
Coordinated Science Lab, University of Illinois at Urbana-Champaign, Urbana, IL\\%
{\tt\small \{solmez2,sa57,jkim684,miehling,basar1,mwest,mehtapg\}@illinois.edu}
}

\begin{document}

\tikzstyle{circ} = [draw,circle,fill = white!20,minimum width = 3pt, inner sep  = 5pt]
\tikzstyle{line} = [draw, -latex']

\maketitle
\thispagestyle{empty}
\pagestyle{empty}
		\begin{abstract}
Evolution of disease in a large population is a function of the top-down policy measures from a centralized planner, as well as the
self-interested decisions (to be socially active) of individual agents
in a large heterogeneous population.  This paper is concerned with
understanding the latter based on a mean-field type optimal
control model.  Specifically, the model is used to investigate the
role of partial information on an agent's decision-making, and study the impact of such
decisions by a large number of agents on the spread of the virus in the population.  The motivation comes from
the presymptomatic and asymptomatic spread of the COVID-19 virus
where an agent unwittingly spreads the virus. We show that even in a setting with fully rational agents, limited information on the viral state can result in an epidemic growth.

	\end{abstract}

\def\b{\beta}
\def\v{\alpha}
\def\f{\bar{\phi}}

\section{Introduction}\label{sec:intro}

Social behavior of individual agents ultimately drives the case count
in an epidemic.  Decisions of an individual agent are mediated by the
top-down policy guidelines, e.g., lockdown mandates 
or social distancing guidelines.  However, these decisions are also a
function of: (i) the agent's risk-reward trade-off based on the agent's health
condition, age, and information
from media on transmissibility and lethality of the virus
(e.g., omicron variant is more transmissible but less deadly than the delta
variant), (ii) aggregate positivity rate in the local population, and (iii) the agent's
assessment of its own epidemiological status.

This paper is a continuation of prior work from our group on modeling
an agent's decision-making in an epidemic based on mean-field game (MFG) formalisms~\cite{olmez2021pre}. The primary focus of the current paper is to investigate the role of agent rationality and information on agent behavior. The motivation comes from the following
thought experiment: On the planet of Vulcan, the agents are both
perfectly rational and perfectly informed.  Because an agent has
perfect information, it immediately knows its epidemiological status.
And because an agent is perfectly rational, it knows to
self-quarantine if infected.  Therefore, the rational choice of a
susceptible (non-infected) agent is to carry on with normal life.  This is
because it can count on rational choices of perfectly informed
others. 

Unfortunately in reality, on the planet of Earth, we continue to experience waves of infection
which are doubtlessly caused by emergence of variants but surely
exacerbated by our own actions which are often less than rational and
certainly almost always based on incomplete information.  A case in
point is the role played by the presymptomatic and asymptomatic spread of the COVID-19 virus~\cite{rivett2020screening,buitrago2020role,bender2021analysis,wei2020presymptomatic}.  

The contributions of this paper are as follows. Apart from the
development of the model, we derive several qualitative insights into an
agent's actions in the midst of an epidemic.  Specifically, these
insights are related to: (i) the risk-reward tradeoff of a susceptible
agent, and (ii) the evolution of belief of an active agent who does
not show symptoms.  The resulting decision-making and its effect on the disease outbreak in the population are discussed at length. These qualitative insights are drawn from quantitative
results presented as part of Props.~\ref{prop:sta_sol}, \ref{prop:partial_threshold}, and \ref{prop:filter}, and as illustrated in Figs.~\ref{fig:risk-reward}-\ref{fig:abar_at_beta05}.  The discussion appears in
\Sec{sec:discuss}.  

Closely related to our work are~\cite{elie2020contact}, where an agent's
decision variable is its rate of contact with others, and~\cite{aurell2020optimal,aurell2021finite}, where an agent strives to
follow a prescribed rate of contact based on government guidelines.
Other MFG-style modeling of epidemics appears in~\cite{hubert2020incentives,lee2020controlling,cho2020mean,doncel2020mean,tembine2020covid}. 
The novelty of our work comes from the inherent partial observability of viral status
and differences in cost structures. These factors are crucial for modeling the asymptomatic spread of an epidemic.

The remainder of the paper is organized as follows. The problem formulation appears in \Sec{sec:prob} where the two main questions of the paper are also introduced.  The answers to these questions are based on the analysis of the
HJB equations derived in~\Sec{sec:optimality_equations}. The
analysis is presented in \Sec{sec:analysis} together with 
qualitative answers to the two main questions.  These answers conclude this paper and spur the development of the MFG formalism. This formalism, which will be studied in detail as part of future work, is presented in \Sec{sec:mfg}. The proofs appear in the Appendix.

\begin{figure*}[!t]
	\centering
	\begin{tabular}{ccc}
		\begin{subfigure}[b]{0.28\textwidth}
		\centering
 		\begin{tikzpicture}[node distance = 2cm,auto]
        \node[circ] (S) {$\Sem$};
        \node[circ, above of = S] (A) {$\Aem$};
        \node[circ, right of = A] (I) {$\Iem$};
        \node[circ, right of = I] (D) {$\Dem$};
        \node[circ, below of = I] (R) {$\Rem$};
        \path[line] (S) -- node[midway,left] {$\eta+ \lambda^{\text{\tiny SA}} \beta_t  U_t$} (A);
        \path[line] (A) -- node[midway,above] {$\lambda^{\text{\tiny AI}}$} (I);
        \path[line] (A) -- node[midway,left] {$\lambda^{\text{\tiny AR}}$} (R);
        \path[line] (I) -- node[midway,above] {$\lambda^{\text{\tiny ID}}$} (D);
        \path[line] (I) -- node[midway,right] {$\lambda^{\text{\tiny IR}}$} (R);
        \node at (-1.7,2.0) {asymptomatic};
        \node at (-1.3,0) {susceptible};
        \node at (3.2,0) {recovered};
        \node at (4.9,2) {dead};
        \end{tikzpicture}
        \caption{}
 		\end{subfigure}
		&$\quad\quad\quad\quad\quad\quad$&
		\begin{subfigure}[b]{0.38\textwidth}
		\centering
		\begin{tabular}{|c|c|c|c|}
		\hline
		\multirow{3}{*}{\textbf{state}} & \textbf{health} & \textbf{alturistic} & \textbf{activity} \\
		& \textbf{cost} & \textbf{cost} & \textbf{reward} \\ \cline{2-4}
		& $c^{\text h}(x)$ & $c^{\text a}(x)$ & $r(x)$ \\
		\hline
		$\Sem$ & $0$ & $0$ & $\v$ \\
		$\Aem$ & $0$ & $1$ & $\v$\\
		$\Iem$ & $1$ & $1$ & $\v$ \\
		\hline
		\end{tabular}
		\caption{}
		\end{subfigure}
	\end{tabular}
\caption{(a) Epidemiological states and the associated transition graph. (b)
  Model for running cost: It is assumed that the economic reward of 
  $\alpha$ is smaller than the altruistic cost of $1$ for an
  asymptomatic ($x=\Aem$)
  or symptomatic  ($x=\Iem$) agent.}
\label{fig:model}
\end{figure*}
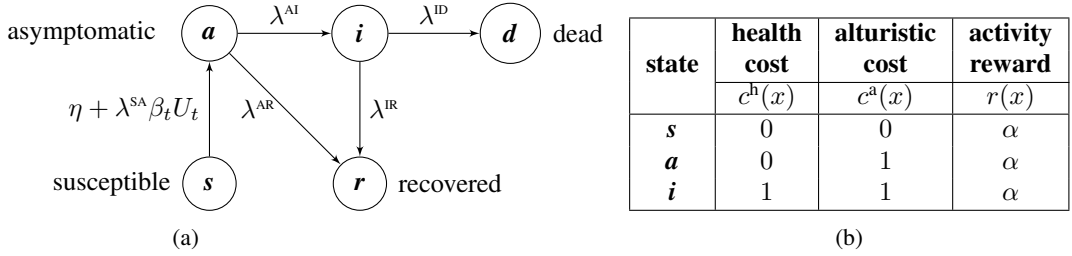

\section{Problem Formulation}\label{sec:prob}

\subsection{Model for a Single Agent}
\newP{Dynamics} The epidemiological state of a single agent is modeled
as a Markov process $X:=\{X_t\in\Xs: t\geq 0\}$ with state-space $\Xs := \{\Sem,\Aem,\Iem,\Rem,\Dem\}$.
Fig.~\ref{fig:model} depicts the transition graph along with a
description of the epidemiological meaning of each of the five
states.  Notably, there are two types of infected states: (i)
asymptomatic state, denoted as $\Aem$; and (ii) symptomatic state, denoted as
$\Iem$.  In either of these states, an agent is 
infectious, i.e., able to infect other agents. The modeling distinction is
that, in partially observed
settings, an asymptomatic agent may not {\em know} its true state but
a symptomatic agent does. Note that unlike our previous
paper~\cite{olmez2021pre}, the transition graph now includes an edge
from $\Aem$ to $\Rem$. This is more realistic, but
significantly complicates modeling and analysis.  Broadly, there are two types of transitions:

\smallskip

\noindent \textbf{1)} On the subset $\{\Aem,\Iem,\Rem,\Dem\}$, the
transition rate depends {\em only} upon the agent attribute $\theta$,
which here represents the age of the agent.  For example, an older
symptomatic agent has a longer expected recovery time (smaller
$\lambda^{\text{\tiny IR}})$ than a younger one.

\smallskip

\begin{figure}[h!]
\includegraphics[width=8cm]{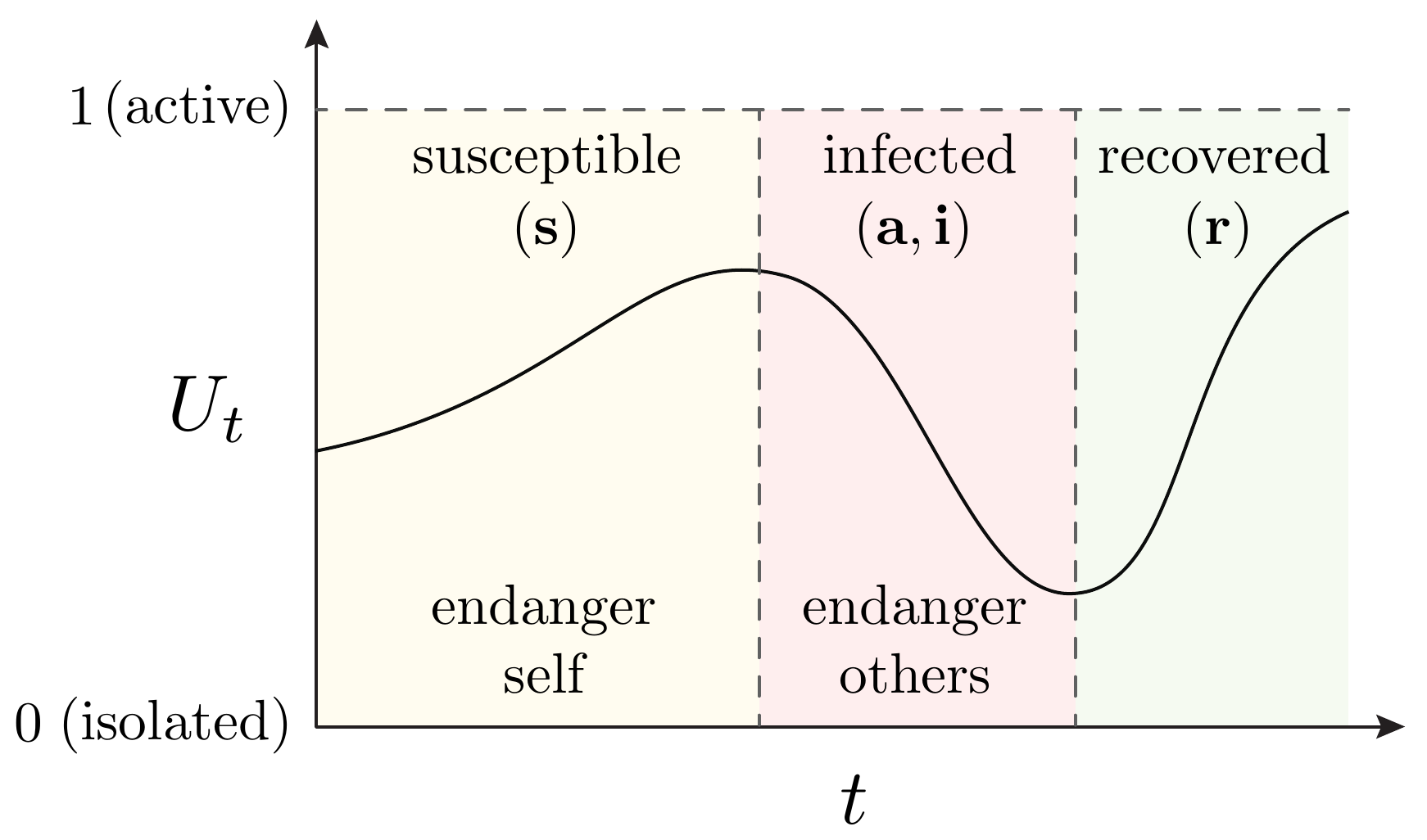}
\caption{A qualitative representation of the control input. Agents choose their activity levels based on their own viral status and the disease level within the population.}
\label{fig:intro}
\end{figure}

\noindent \textbf{2)} The transition from $\Sem \to \Aem$ depends upon three factors:
  (i) the intrinsic infectivity of the virus, (ii) the agent behavior
  (level of social activity), and (iii) the behavior of the
  infected agents in the population. The following equation is used
  to model the effect of these three factors:
 \[
 \text{rate}[\Sem \to \Aem] = \eta + \lambda^{\text{\tiny SA}}\ \b_t\,  U_t 
 \]
 where $\eta>0$ is a small baseline rate (useful for regularizing the
 problem) and $\lambda^{\text{\tiny SA}}$ is the virus transmissibility
 parameter.  
 The process $U$ is referred to as the agent's control input, and
 models the agent's activity {level}:  With
 $U_t=0$ (resp., $U_t=1$) the agent is isolated
 (resp., normally active) at time $t$ (Fig.~\ref{fig:intro}).  The process $\beta$ models
 the average activity level of the infected agents
 (Eq.~\eqref{eq:mf-vars}).

The two main questions driving this work are as follows: (i) How does
a single agent choose its control input $U$? and (ii) How does that
choice (made by agents) affect the evolution of disease in a large heterogeneous
population?
To answer these questions, we adopt an optimal control framework.

\medskip

\noindent \textbf{Optimal control objective:}  In the following, $\b$
is a given deterministic process. The control objective for a single
agent is to choose its activity $U$ to minimize 
\[
{\sf J}(U;\b)  = {\sf E} \left( \int_0^T
  e^{-\gamma t} c(X_t,U_t) \ud t + e^{-\gamma T} \phi(X_T) \right)
\]
where $\gamma > 0$ is a discounting rate (assumed to be much smaller than the transition rates),
$T=T(\omega)=\inf\{t>0: X_t(\omega) \in \{\Dem,\Rem\}\}$ is the random
stopping time when the agent either recovers ($X_T=\Rem$) or the agent
dies ($X_T=\Dem$); by convention, $\inf \emptyset = \infty$. The cost
function is of the following form:
%
\begin{align*}
c(x,u) &= c^{\text{h}} (x) +\left( c^{\text{a}} (x) - r(x) \right) u
\end{align*}
where $c^{\text{h}}(\cdot)$, $c^{\text{a}}(\cdot)$, and $r(\cdot)$
model the health cost, altruistic cost, and economic reward,
respectively (Fig.~\ref{fig:model}~(b)). The altruistic cost models an agent's desire to help the greater good and has the effect of inducing an infected agent to isolate. Of course, an agent's
choice is mediated by its information, the structure of which is described next.



 


\medskip

\noindent \textbf{Information structure:} There are two settings of
the problem: (i) the fully observed case, and (ii) the partially
observed case. 
In the partially observed setting, the observation process 
$Y:=\{Y_t\in\{0,1,2\}:t\geq 0\}$ is defined as
\[
Y_t = 1_{[X_t=\Iem]} + 2 \cdot 1_{[X_t=\Dem]}
\]
Therefore, an infected agent who is also asymptomatic (in state
$\Aem$) does not observe its epidemiological status until it begins to
show symptoms (in state $\Iem$).  

\subsection{Model for the Mean-Field}
\label{sec:pomfg_model}

To fully specify the problem, we need to define the model for the agent's attribute $\theta$ and the process $\b$, henceforth referred to as the
mean-field process.   
The probability mass function of the attribute $\theta$ is denoted by
${\sf p}(\cdot)$. Denote by $\rho_t(x,u;\theta)$ the joint
distribution of the state-action pair $(X_t,U_t)$ at time $t$,
conditioned on the attribute $\theta$.  At time $t$, the average
activity level of infected agents is 
\begin{equation}\label{eq:mf-vars}
\b_t := \sum_{\theta}  {\sf p} (\theta) \sum_{x\in\{\Aem,\Iem\}} 
\int_0^1 u 
              \rho_t(x,u;\theta) \ud u \\
\end{equation}
which we define to be the mean-field process.


\subsection{Function Spaces}\label{sec:fnspaces}

The filtration of the Markov
process $X$ is denoted by ${\cal F}:=\{{\cal F}_t:t \geq 0\}$ where
${\cal F}_t:=\sigma(X_t)$.  The filtration of the observation process
$Y$ is denoted by ${\cal Y} := \{ {\cal Y}_t:t\geq 0\}$ where ${\cal
  Y}_t:=\sigma(\{Y_s:0\leq s\leq t\})$.  In the two settings of the
problem, the spaces of admissible control inputs, denoted by ${\cal U}$ in each case, are as
follows:
\begin{align*}
\text{(fully obsvd.)} \qquad {\cal U} &= L^2_{\cal F}([0,\infty) ;[0,1])\\
\text{(part. obsvd.)} \qquad {\cal U} &= L^2_{\cal Y}([0,\infty) ;[0,1])\
\end{align*}
i.e., an admissible control input $U$ is a $[0,1]$-valued stochastic
process adapted to ${\cal F}$ in the fully observed case, and
adapted to ${\cal Y}$ in the partially observed case. The use of the
common notation ${\cal U}$ should not cause any confusion because the
two cases are treated in separate subsections. The process $\beta$ is assumed to be deterministic. The function space
for $\beta$ is ${\cal M}=L^2([0,\infty);[0,1])$.

\section{Optimality equations for a single agent}
\label{sec:optimality_equations}

\subsection{Fully Observed Case}
\label{sec:fully}

For each $x\in\Xs$ and $t\geq 0$, the value function is
\begin{align}
& v_t(x) :=  \min_{U\in L^2_{\cal F}} \nonumber \\ & {\sf E} \left( \left. 
         \int_t^{T} e^{-\gamma(s-t)}c(X_s,U_s) \diff s
        +e^{-\gamma T} \phi(X_T) \right| X_t= x\right)
\label{eq:val_fully}
\end{align}
where at the two terminal states $x=\Dem$ and $x=\Rem$, the value function is given by known constants $v_t(\Dem)=\f(\Dem)$ and
$v_t(\Rem)=\f(\Rem)$.



For the fully observed problem, complete characterizations of the
value function and the optimal control are described in the following
proposition.  This result is a minor extension of a similar result
appearing in our prior paper~\cite{olmez2021pre} (for
$\lambda^{\text{\tiny AR}}$ and $\eta=0$) and therefore its proof is
omitted.  



\medskip

\begin{proposition}
Suppose $\alpha<1$. 
For $x \in \{\Aem,\Iem,\Rem,\Dem\}$, the value function $v_t(x)$ and
the optimal control law $\psi_t(x)$ are
stationary as tabulated in Table~\ref{tab:val_function}. For state
$\Sem$, the value function $v_t (\Sem)$ solves the HJB equation
\begin{align*}
-\frac{\ud v_t}{\ud t} (\Sem)& + (\gamma+\eta) v_t (\Sem) \\
&= \eta \bar{\phi} (\Aem) + \min_{u \in [0,1]} \left( \be \left( \bar{\phi} (\Aem) - v_t (\Sem) \right) - \v \right) u
\end{align*}
\label{prop:prop1}
\end{proposition}

{\renewcommand{\arraystretch}{2.5}
\begin{table}[h!]
\centering
\begin{tabular}{|c|c|c|}
\hline
\textbf{state} $x$ & \textbf{value function} $v_t(x) $& \textbf{optimal
                                                   control}\\
\hline
$\Aem$ &{\large $\frac{c^h(\Aem)+\la\bar{\phi} (\Iem) +\lar\bar{\phi}
         (\Rem)}{\gamma+\la+\lar}$} $=:\bar{\phi}(\Aem)$ &
                                                           $U_t^{\text{opt}} = \psi_t(\Aem) =0$\\
$\Iem$ &{\large $\frac{c^h(\Iem) + \lr \bar{\phi} (\Rem) + \ld
         \bar{\phi} (\Dem)}{\gamma+\lr+\ld}$} $=:\bar{\phi}(\Iem)$ & $U_t^{\text{opt}}= \psi_t(\Iem) =0$ \\
$\Rem$ & $\bar{\phi} (\Rem)$ (given) & - \\
$\Dem$ & $\bar{\phi} (\Dem)$ (given) & -  \\
\hline
\end{tabular}
\caption{Fully observed case.}
\label{tab:val_function}
\end{table}
{\renewcommand{\arraystretch}{1}

\begin{remark}
The optimal control for both asymptomatic and symptomatic agents is
zero.  This is entirely because $\alpha<1$.  Recall that $\alpha$ is a model
for economic reward per unit time.  The assumption $\alpha<1$ means
that the economic reward ($\alpha$) is outweighed by the altruistic cost ($1$).
\end{remark}

\subsection{Partially Observed Case}

The partially observed problem is first converted to a fully observed
one by introducing the belief state, as in~\cite{saldi2019partially,saldi2019po}, which at time $t$ is denoted by
\[
\pi_t := \begin{bmatrix} \pi_t(\Sem) & \pi_t(\Aem) &
  \pi_t(\Iem)  & \pi_t(\Rem) & \pi_t(\Dem) \end{bmatrix}
\]
where $\pi_t(x):= {\sf P} ([X_t = x] \mid {\cal Y}_t)$ for $x \in \Xs$.
Since the events $[X_t = \Iem]$ and $[X_t = \Dem]$ are both contained
in ${\cal
  Y}_t$, $\pi_t$ is not an arbitrary element of the
probability simplex in $\mathbb{R}^5$.  Let ${\cal P}^1$ denote the
set of pmfs on $\{\Sem,\Aem,\Rem\}$ and let ${\cal P}^2
= \{\delta_\Iem,\delta_\Rem,\delta_\Dem\}$.  Then, the state-space for
the belief is ${\cal P}^1 \cup {\cal P}^2$. For $t\geq 0$ and 
$\mu\in {\cal P}^1 \cup {\cal P}^2$, the value function is
\begin{align*}
& v_t(\mu)  :=\min_{U\in L^2_{\cal Y}}\\
& {\sf E} \left(
           \int_t^Te^{-\gamma(s-t)} c(X_s,U_s;\v) {\rm
           d} s + e^{-\gamma T}\phi(X_T) \mid \pi_t = \mu  \right)
\label{eq:val_partially}
\end{align*}
There are two cases to consider: (i) when $\mu\in{\cal P}^2$, and (ii)
when $\mu\in{\cal P}^1 \setminus \{ \delta_\Rem \}$.  In the first case, when $\mu\in{\cal P}^2$,
the problem reduces to the fully-observed setting, and the value
function is given by
\[
v_t(\delta_\Dem) = \f(\Dem),\;\;v_t(\delta_\Rem) = \f(\Rem),\;\;v_t(\delta_\Iem) = \f(\Iem)
\] 
The optimal control for the agent in the symptomatic state
($\pi_t=\delta_{\Iem}$) is $U_t^{\text{opt}} = \psi_t (\delta_{\Iem})
= 0$.  

For the second case, when $\mu\in{\cal P}^1 \setminus \{ \delta_\Rem \}$, a nonlinear filter is
used to obtain the evolution of the belief.  For this purpose,
consider first the random variable $\tau=
\tau (\omega) = \inf \{ t>0 : X_t(\omega) = \Iem \}$.  Now, $\tau$ is
a ${\cal Y}_t$-stopping time and 
\[
\pi_t= \begin{bmatrix} \pi_t(\Sem) & \pi_t(\Aem) &
  0 & \pi_t(\Rem) & 0 \end{bmatrix} \quad \text{for}\;\;t<\tau
\]
Let $S_t:= \pi_t(\Sem)$, $A_t:=\pi_t(\Aem)$ and $R_t:= \pi_t(\Rem)$ for $t<\tau$.  Then the
stochastic process $\{(S_t,A_t,R_t)\in[0,1]^3: S_t+A_t+R_t = 1, 0\leq t<\tau\}$ evolves according to the nonlinear
filter
\begin{subequations}
\begin{align}
\frac{\ud S_t}{\ud t} &= \left( -\be U_t - \eta + A_t \la \right) S_t\\
\frac{\ud A_t}{\ud t} &= \left( \be U_t + \eta \right) S_t + A_t \left( -\la - \lar + \la A_t \right)\label{eq:dA_dt}\\
\frac{\ud R_t}{\ud t} &= A_t (\lar+\la R_t)
\end{align}
\label{eq:filter}
\end{subequations}
which is derived from the general form detailed by~\cite{confortola2013filtering}.

\begin{figure}[h!]
\centering
\includegraphics[width=7cm]{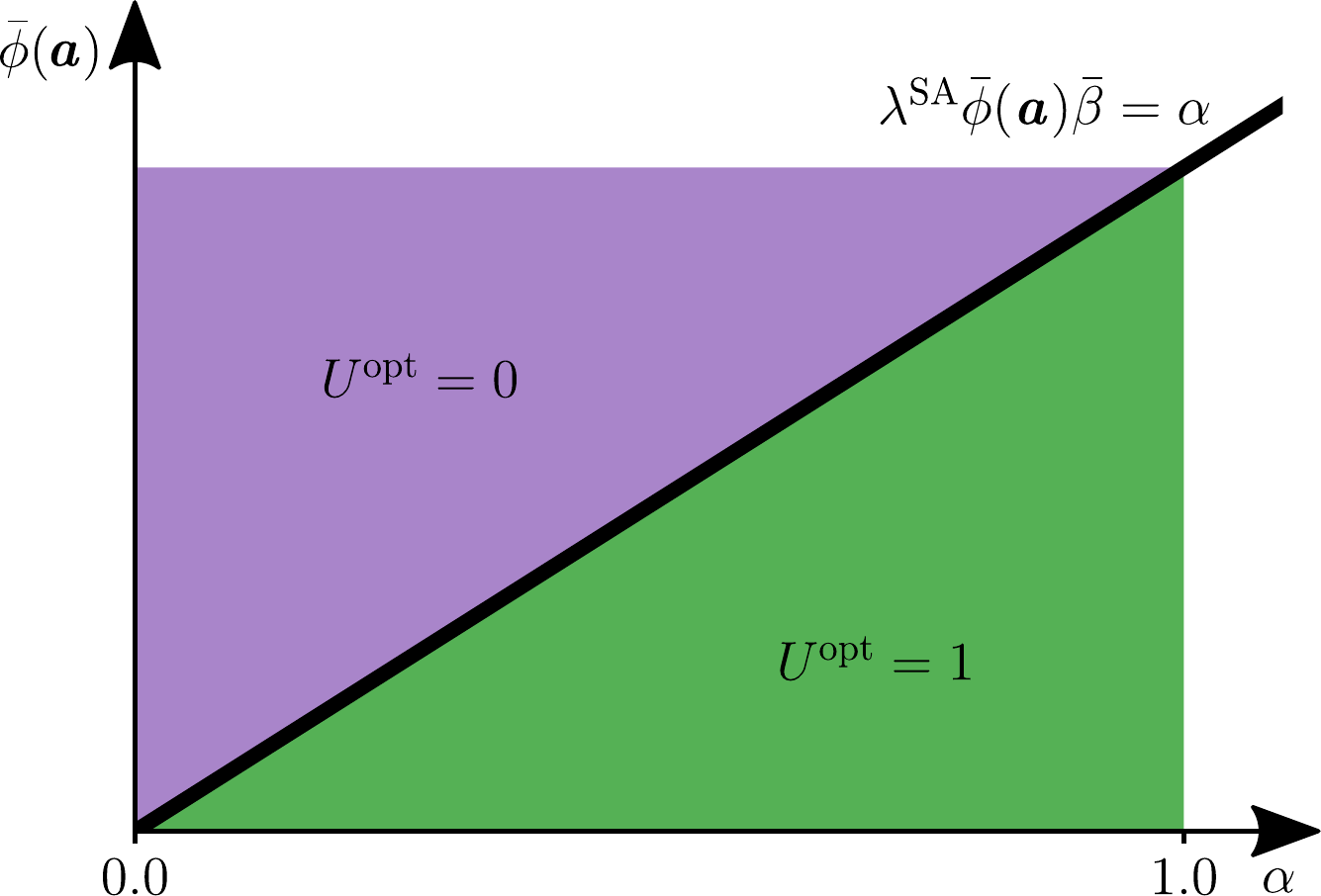}
\caption{A susceptible agent's risk-reward diagram. The line
  $\lambda^{\tiny \text{SA}} \bar{\phi}(\Aem) \bar{\beta} = \alpha$ is
where the risk equals the reward.}
\label{fig:risk-reward}
\end{figure}


We identify ${\cal P}^1$ with the domain $D:= \{(s,a) \in [0,1]^2: s+a \leq 1\}$ where $s$ is the value of $S_t$ and $a$ is the value of $A_t$.  For an arbitrary
element $\mu = [s,a, 0,1-s-a, 0]$ in ${\cal P}^1$, we denote the
value function with respect to the $(s,a)$-coordinates as 
\[
\phi_t(s,a) := v_t(\mu), \quad (s,a) \in D,\quad t\geq 0
\]
The following proposition provides the HJB equation whose derivation
appears in Appendix~\ref{appdx:derivations}.  

\medskip

\begin{proposition} The value function $\phi_t(s,a)$ solves the HJB equation
\begin{align}
\left( -\frac{\partial \phi_t}{\partial t} + \gamma \phi_t\right)(s,a) = ( {\cal L}\phi_t)(s,a)
+ \min_{u \in [0,1]} {\cal M}_t(s,a) u
\label{eq:HJB_POMFG}
\end{align}
together with boundary conditions
\begin{subequations}
\begin{align}
&\phi_t(s,a)-\phi_t(s,1-s-a) = (2a+s-1) \left( \bar{\phi} (\Aem) -
  \bar{\phi} (\Rem) \right) \\
&\phi_t(0,a) = a \bar{\phi} (\Aem) + (1-a) \bar{\phi} (\Rem) 
\end{align}
\label{eq:bcs}
\end{subequations}
where the formulae for the linear operators ${\cal L}$ and ${\cal M}_t$
appear in Appendix~\ref{appdx:derivations}.    
\end{proposition}

\medskip

In the
remainder of this paper, it is assumed that a unique solution for the
HJB equation exists and yields a
well-posed optimal control law, denoted as $\psi_t(s,a)$ for
$(s,a)\in D$.


\section{Analysis of the decision-making by an agent} \label{sec:analysis}


In this section, an analysis of the HJB and the nonlinear filter
equations is given.  The analysis is used to develop insights into the decisions
of a single agent in an epidemic (given $\beta$).  To aid the
analysis, it is assumed that $\beta_t = \bar{\beta}$ for all $t\geq 0$
(i.e. stationary).  The main
results are described in the form of three propositions (Props.~\ref{prop:sta_sol}-\ref{prop:filter}) and
illustrated with accompanying figures (Fig.~\ref{fig:risk-reward}-~\ref{fig:abar_at_beta05}).  The section
concludes with a qualitative 
discussion of an agent's decision-making in~\Sec{sec:discuss}.

\subsection{Risk-Reward Tradeoff for a Fully Observed Susceptible Agent}

\begin{proposition}[Stationary Solution] Suppose $\v<1$,
  $\gamma>0$, $\f(\Iem)>0$, $\b_t \equiv
  \bar{\b}$.  Set $\bar{\beta}^{\text{\tiny crit}} := \frac{\v}{\lambda^{\text{\tiny
                                SA}}\f(\Aem)} \left(
                            1+\frac{\eta}{\gamma} \right)$.  Then the optimal control for a
  susceptible agent is stationary and described by the following cases:
\begin{enumerate}
			\item If $
                          \bar{\beta}<\bar{\beta}^{\text{\tiny crit}}$, then
			the  optimal control $U_t^{\text{opt}}=\psi_t(\Sem)=1$ and the
			optimal value $v_t(\Sem) =
                        \frac{(\lambda^{\text{\tiny SA}} \bar{\beta}+\eta)\bar{\phi} (\Aem)- \v}{\lambda^{\text{\tiny SA}}\bar{\b} + \gamma + \eta}$.  
			\item If $\bar{\b}\geq
                         \bar{\beta}^{\text{\tiny crit}}$, then
			the optimal control $U_t^{\text{opt}}=\psi_t(\Sem)=0$ and the optimal value $v_t(\Sem) = \frac{\eta}{\gamma+\eta} \bar{\phi} (\Aem)$.
\end{enumerate}

\label{prop:sta_sol}
\end{proposition}

\medskip

\begin{remark}
A susceptible agent's decision is best understood using a risk-reward
diagram depicted in Fig.~\ref{fig:risk-reward}.  The horizontal axis of this
diagram is the economic reward per unit time ($\alpha$).  The vertical axis is the
potential health cost (value function $\bar{\phi} (\Aem)$).  The
diagram helps show that the optimal decision for an agent is to be
active if the reward is greater than the risk.  For any given
$\alpha$ and $\bar{\phi} (\Aem)$, the risk-reward analysis reveals
a critical threshold $\bar{\beta}^{\text{\tiny
    crit}}$ above which the agent ceases to be active.  It is noted
that the critical threshold scales inversely with the product
$\lambda^{\text{\tiny SA}}\f(\Aem)$.  This is useful in several ways
for analysis:
\begin{enumerate}
\item An older agent will have a greater potential health cost
  $\f(\Aem)$ and therefore a smaller critical threshold $\bar{\beta}^{\text{\tiny
    crit}}$.  
\item A more contagious (high $\ls$) variant may still have a larger critical threshold, if it is significantly less severe (low $\f(\Aem)$). For instance, the omicron variant is two times more transmissible ($2\times$
  greater $\lambda^{\text{\tiny SA}}$)  than the delta
  variant~\cite{nyt2022where}. However, the potential health cost $\f(\Aem)$ of
  the omicron is also smaller (because it is less lethal).  
\end{enumerate}
\end{remark}


\begin{figure}[b]
\centering
\includegraphics[width=6cm]{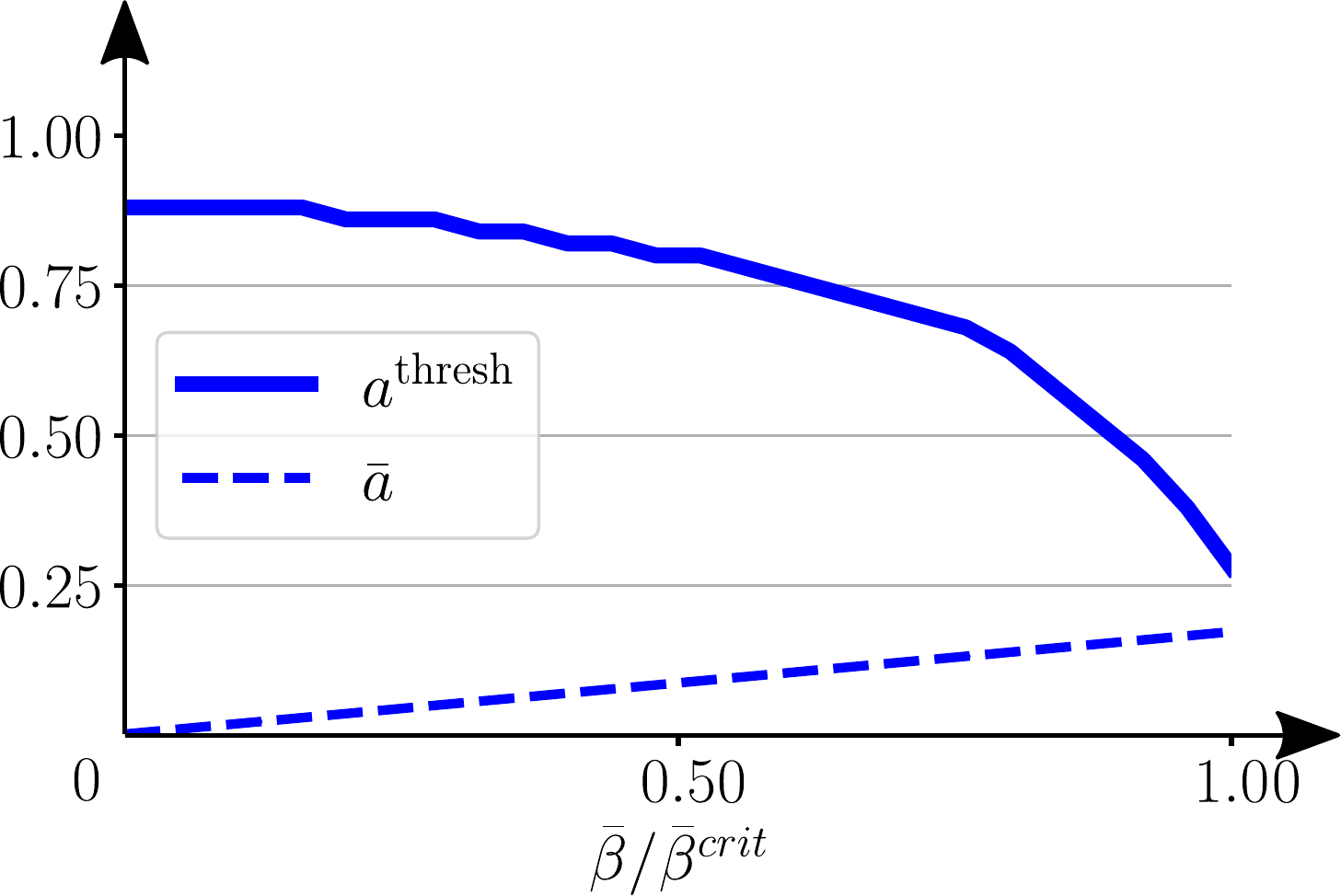}
\caption{Plot of $\at$ and $\bar{a}$ as a function of $\bar{\b} / \bar{\b}^{\text{crit}}$.}
\label{fig:at_abar}
\end{figure}

\begin{figure*}[t]
     \centering
     \begin{subfigure}[b]{0.32\textwidth}
         \centering
         \includegraphics[width=\textwidth]{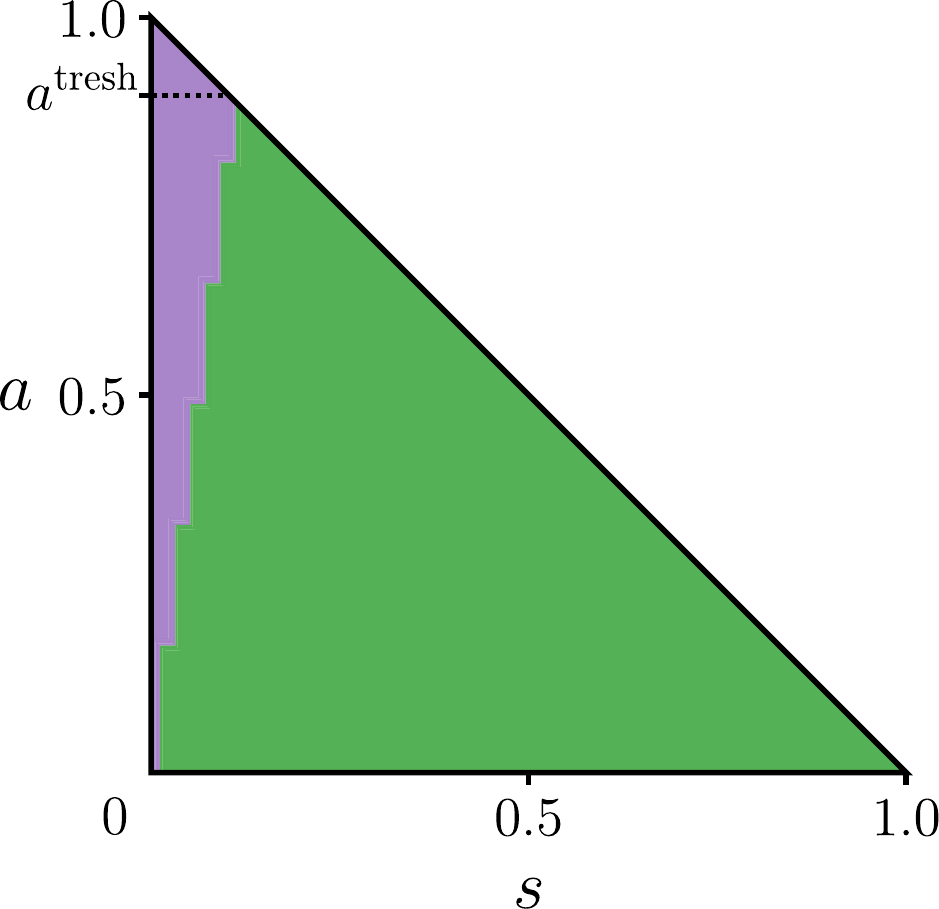}
     \end{subfigure}
     \begin{subfigure}[b]{0.32\textwidth}
         \centering
         \includegraphics[width=\textwidth]{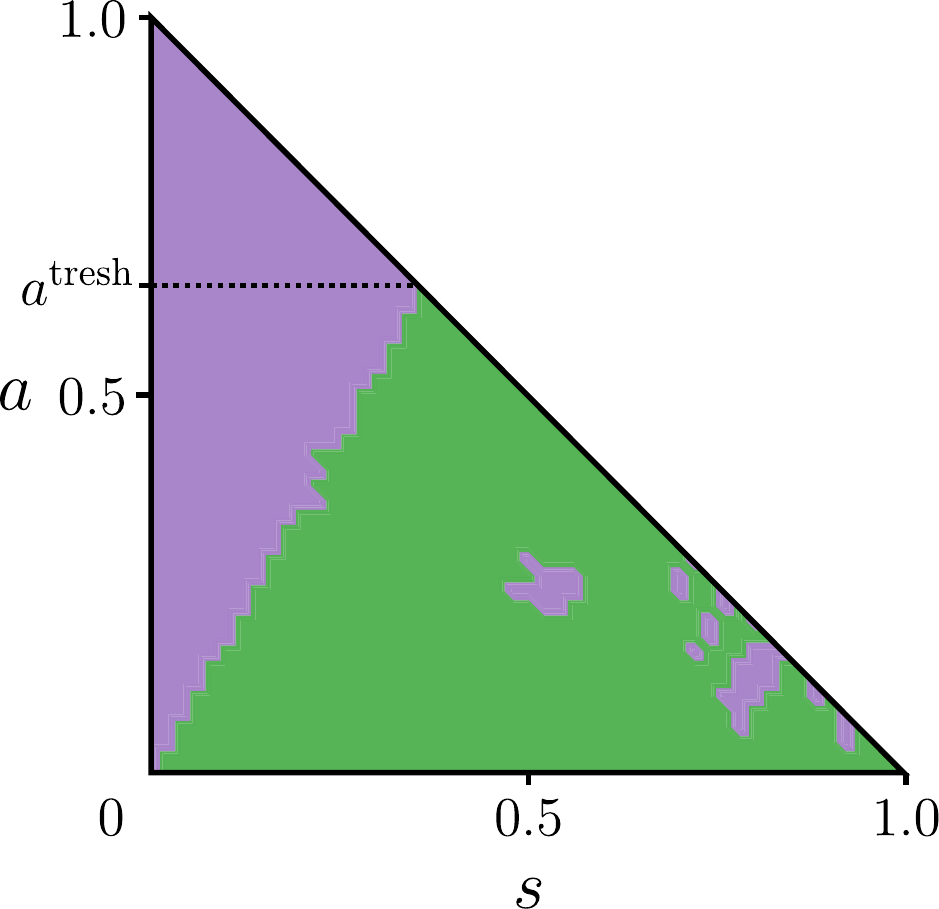}
     \end{subfigure}
     \begin{subfigure}[b]{0.32\textwidth}
         \centering
         \includegraphics[width=\textwidth]{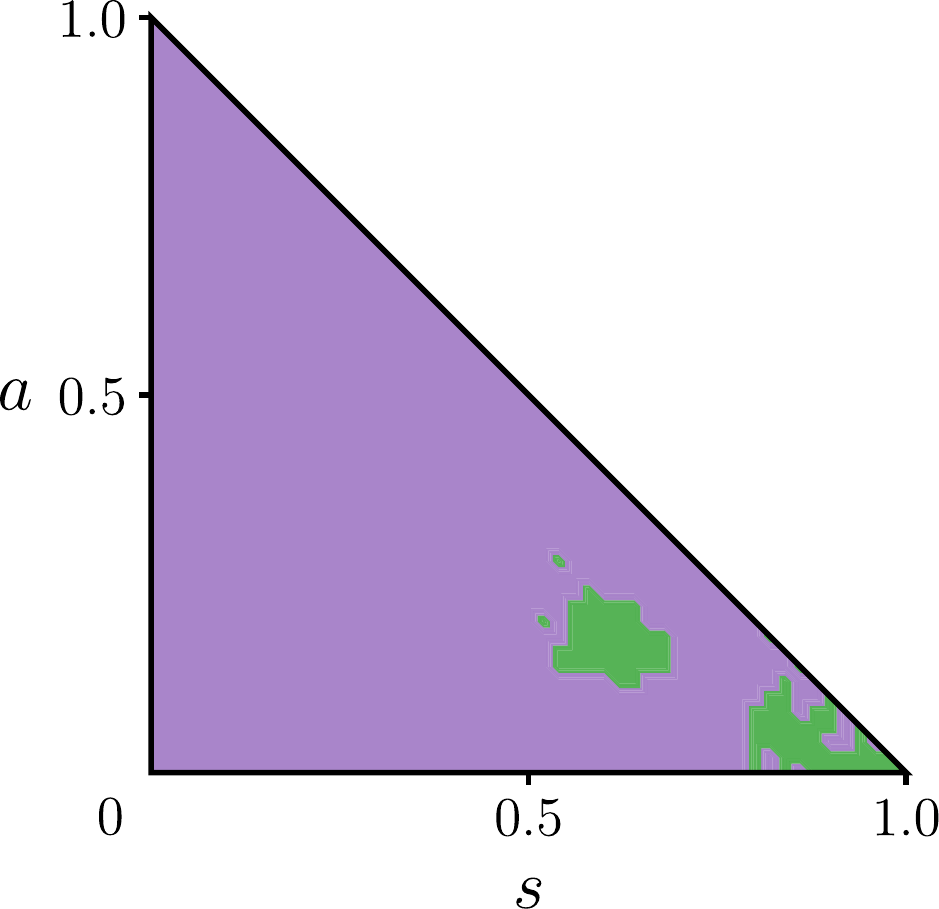}
     \end{subfigure}
        \caption{Stationary optimal control policies for the partially
          observed case; 
          $\alpha = 0.9$ and $\bar{\beta} /
          \beta^{\text{crit}} = 0.1 , 0.8, 1.1$ for the three plots. Purple region is
          $U^{\text{opt}} = 0$ and the green
          region is $U^{\text{opt}} = 1$. The islands are
          numerical artifacts.}
        \label{fig:u_sta_beta}
\end{figure*}

\subsection{Threshold Optimal Policies for Partially Observed Case}

In our prior paper~\cite{olmez2021pre}, we derived the following result
for the solutions of the HJB equation in a special case:  

\begin{proposition}\label{prop:partial_threshold}
Suppose $\b_t \equiv
  \bar{\b}$, $S_0 + A_0 =1$, and $\lar = 0$.  Then, for $t<\tau$, $S_t + A_t =1$, and
\begin{enumerate}
\item If $\bar{\b}\geq \bar{\b}^{\text{crit}}$,
  then the optimal control law ${\bar{\psi}}(\cdot,a) = 0$.
\item For each fixed  $\bar{\beta}<\frac{\v}{\lambda^{\text{\tiny
        SA}}\f(\Iem)}<
  \bar{\b}^{\text{crit}}$ there exists a
  $\underline{\lambda}^{\text{\tiny AI}} = \underline{\lambda}^{\text{\tiny AI}}
  (\bar{\b})$ such that for all $\lambda^{\text{\tiny AI}} >
  \underline{\lambda}^{\text{\tiny AI}}$, the optimal control law is of
  threshold type:
\begin{align*}
\bar{\psi}(\cdot,a)= \begin{cases}
1 & \text{if} \;\; 0\leq a< {a}^{\text{thresh}} \\
0 & \text{if} \;\; {a}^{\text{thresh}} < a \leq 1
\end{cases}
\end{align*}
where the threshold $ {a}^{\text{thresh}}\in(0,1)$.  
Furthermore $\ls \bar{\beta} / \la <{a}^{\text{thresh}} $. The function $\underline{\lambda}^{\text{\tiny AI}}
  (\bar{\b})$ is monotonic in its argument and $\lim_{\bar{\b}\downarrow 0}
    \underline{\lambda}^{\text{\tiny AI}}
  (\bar{\b})=0$.  
  
\end{enumerate}
\end{proposition}

\medskip

With $\lar > 0$, an analytical treatment has so far not been
possible.  However, numerical evidence suggests that the
policies are also of threshold type even with $\lar > 0$.  
Three such numerically obtained stationary policies are depicted in
Fig.~\ref{fig:u_sta_beta}.  The policies were computed using the method of lines
numerical algorithm described in Appendix~\ref{sec:numol}.  From these numerically
computed policies, a threshold of ${a}^{\text{thresh}}$ is identified such
that an agent becomes inactive for $A_t > {a}^{\text{thresh}}$ (see
Fig.~\ref{fig:u_sta_beta}).  By identifying the thresholds for different choices of
$\bar{\beta} < \bar{\beta}^{\text{\tiny
    crit}}$, a plot of ${a}^{\text{thresh}}$ as a function of $\frac{\bar{\beta}}{\bar{\beta}^{\text{\tiny
    crit}}}$ is obtained as depicted in Fig.~\ref{fig:at_abar}.  As $\bar{\beta}
\uparrow \bar{\beta}^{\text{\tiny
    crit}}$, the threshold  ${a}^{\text{thresh}} \downarrow 0$ and the agent
ceases to be active.



\subsection{Belief of an Active Agent who Shows No Symptoms}

The analysis thus far has revealed two important insights captured by
the critical value for $\bar{\beta}$ and the threshold for $a$: 
\begin{enumerate}
\item A susceptible agent in the fully observed case is active if $\bar{\beta}$ is
  small enough ($\bar{\beta}<\bar{\beta}^{\text{\tiny
    crit}}$). 
\item In the absence of symptoms, a partially observed agent is active
  if its belief $A_t$ is small enough ($A_t < {a}^{\text{thresh}}$).  
\end{enumerate}

Therefore, to determine the actions of an agent in the partially observed case, it
becomes important to understand the evolution of the process
$\{A_t: 0\leq t< \tau\}$.  This is the subject of
the following proposition.   


\medskip

\begin{proposition}
Set $\bar{a} :=
\frac{\beb+\eta}{\la}$.  Suppose $t<\tau$, and
\begin{align*}
S_0 & = 1 \qquad & \text{(agent starts
out as susceptible)} \\
U_s & = 1 \quad \forall \; s<t \quad & \text{(agent has been active
                                       upto time} \; t\text{)} 
\end{align*}
Then, 
\[
A_t < \bar{a}, \qquad 0\leq t < \tau
\]   
Moreover if $\lar = 0$, then $
A_t \uparrow \bar{a}$ as $t\to \infty$. 
\label{prop:filter}
\end{proposition}

\medskip

A numerical illustration of the result of this proposition appears as
part of Fig.~\ref{fig:filter}.  The most important point about this plot is that
the agent's belief that it is asymptomatic scales as
$\eta+\mathcal{O}(\bar{\beta})$.


\subsection{Qualitative Discussion}\label{sec:discuss}

We are now ready to provide a qualitative answer to the two primary questions we had raised: (i) How does
a single agent choose its control input $U$?; and (ii) How does that
choice (made individually and independently by many agents) affect the evolution of disease in a large heterogeneous
population?

\medskip

The answer to the first question is as follows:
\begin{enumerate}
\item In the fully observed setting, a susceptible agent is active
  if $\bar{\beta}<\bar{\beta}^{\text{\tiny
    crit}}$.  In this setting, an asymptomatic agent isolates (Table~\ref{tab:val_function}).
\item In the partially observed setting, an agent's decision is based
  on its belief $(S_t,A_t,R_t)$.  An agent who starts out as
  susceptible and active has a belief $A_t <
  \bar{a} = \eta + \mathcal{O}(\bar{\beta})$ (Fig.~\ref{fig:filter}). 
\item For small values of $\bar{\beta} < \bar{\beta}^{\text{\tiny
    crit}}$, the quantity $\bar{a} < {a}^{\text{thresh}}$
(Fig.~\ref{fig:at_abar}).  Therefore, an agent's optimal decision is to continue
to be active. This is illustrated in Fig.~\ref{fig:abar_at_beta05} for $\bar{\b} / \bar{\b}^{\text{crit}} = 0.5$.
\end{enumerate}

\begin{figure}[t]
\centering
\includegraphics[width=8cm]{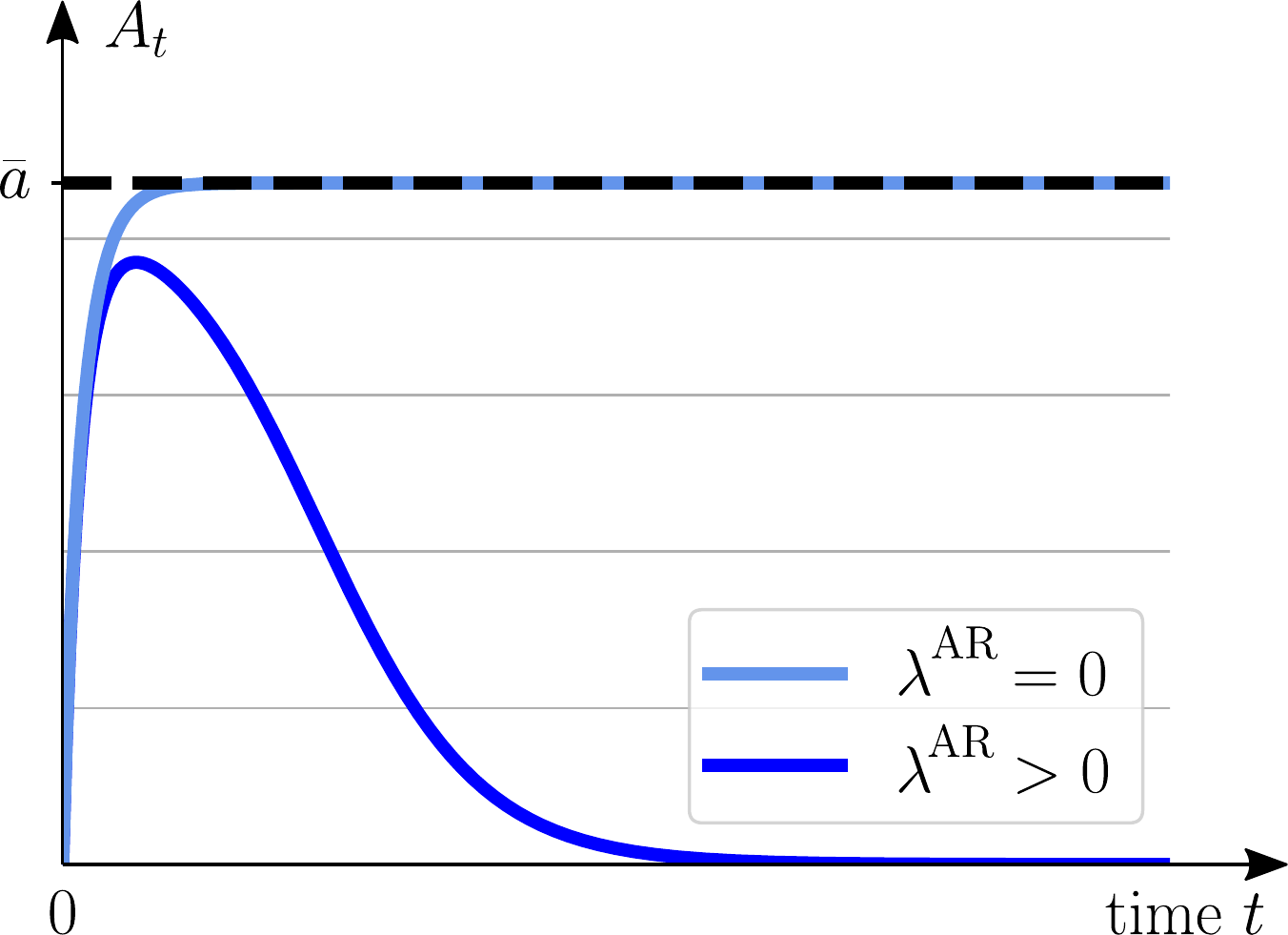}
\caption{Evolution of $A_t$ for an active agent who starts out as susceptible.}
\label{fig:filter}
\end{figure}

\begin{figure}[t]
\centering
\includegraphics[width=6cm]{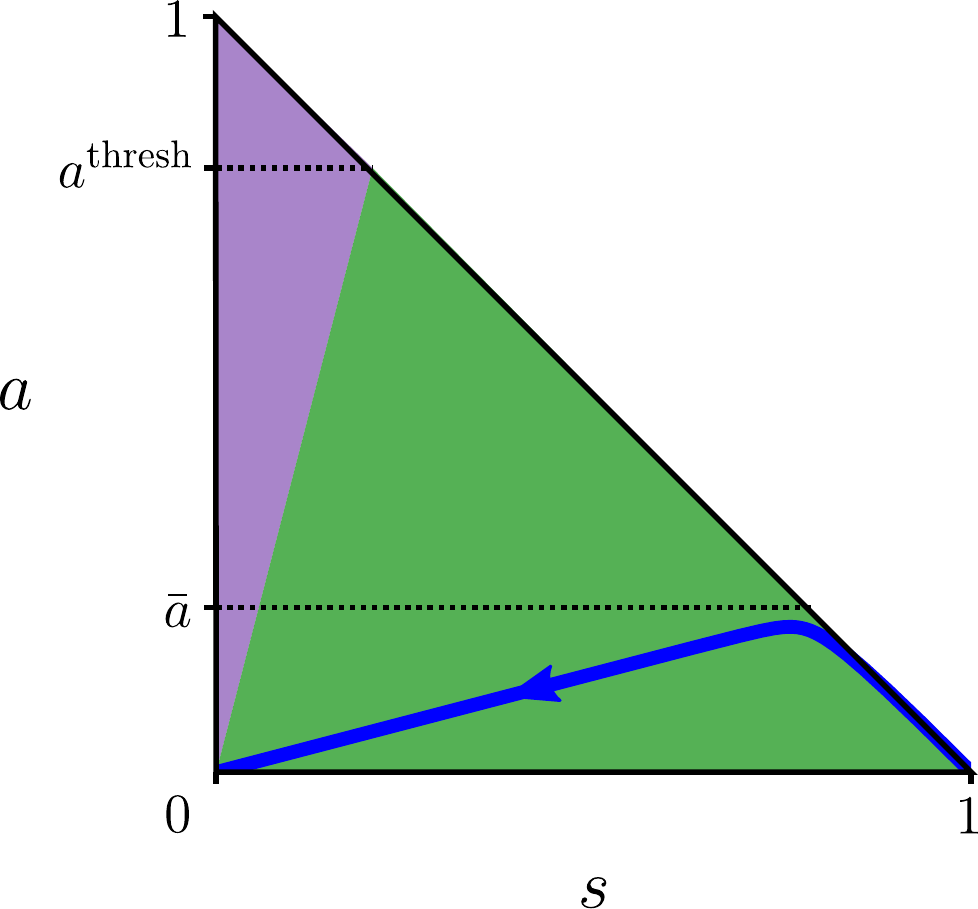}
\caption{The trajectory of Fig.~\ref{fig:filter} (for $\lar>0$)
  depicted on the simplex. An agent’s optimal decision is to remain active, unless it starts to show symptoms.}
\label{fig:abar_at_beta05}
\end{figure}

The answer to the second question is as follows:
\begin{enumerate}
\item Even though an individual agent's conditional probability $A_t$
  is small at any time $t<\tau$, the law of large numbers (LLN)
  dictates that in a large population
  with $N$ agents, ${\sf E}(A_t)$ fraction of agents are asymptomatic
  with high probability.  
\item The basic reproduction
number $R_o$ (pronounced $R$ nought) is defined as $\frac{T^{\text{r}}}{T^{\text{c}}}$ where $T^{\text{r}}$ is the
average time until removal (in our settings, either showing symptoms
$x=\Iem$ or recovery $x=\Rem$) and 
$T^{\text{c}}$ is the average time between infectious contacts.  For
our model, this ratio is $\frac{\ls}{\la+\lar}=R_o$.
\item If $R_o>1$ then for small values of $\bar{\beta}$, optimal
  actions of individual agents (to remain active) still causes an
  epidemic.  
\end{enumerate}
The above is an example of rational irrationality where rational
choices of individual agents lead to irrational outcomes for the
population~\cite{cassidy09rational,Yin_synch_game}.  It is notable that the irrationality here is entirely
a consequence of the lack of information (with full information,
asymptomatic agents isolate).  

A limitation of these arguments is that the analysis is rather
qualitative and moreover relies on the assumption of a given stationarity
$\beta_t\equiv\bar{\beta}$.  Indeed, this is not necessarily true in the case of
an epidemic: The process $\beta$ is non-stationary and also a result
of the individual actions of many agents.  This motivates the need for
a more refined analysis based on a mean-field game model as described next.

\section{Mean-field game model}\label{sec:mfg}

In this section, we consider a population of heterogeneous
agents.  As already noted, $\rho_t(x,u;\theta)$ is the joint
distribution of the state-action pair $(X_t,U_t)$ at time $t$,
conditioned on the attribute $\theta$. The marginal pmf
$\rho_t(x;\theta)$ evolves according to the mean-field equations:
\begin{subequations}
\begin{align}
\frac{\ud \rho_t}{\ud t} (\Sem;\theta) & = -\lambda^{\text{\tiny SA}}
                                        \b_t \int_0^1 u 
                                         \rho_t(\Sem,u;\theta) \ud u -\eta \rho_t (\Sem;\theta)
                                         \label{eq:rho-s} \\
\frac{\ud \rho_t}{\ud t} (x;\theta) & = ({\cal
                                      A}^\dagger\rho_t)(x;\theta),\quad
                                      x\in\{\Iem,\Aem,\Rem,\Dem\} \label{eq:rho-rest}
\end{align}
\label{eq:rho}
\end{subequations}
from a given initial condition $\rho_0(x;\theta)$; 
${\cal A}^\dagger$ is the adjoint of the generator ${\cal A}$ of the
Markov process $X$.  

Each agent in the population chooses its control according to the optimal
control law $\psi_t(\cdot;\theta)$ obtained from solving the HJB
equation (the notation 
denotes the dependence on the attribute $\theta$).  In the two cases,  
\begin{align*}
\text{(fully obsvd.)} \qquad U_t & = \psi_t(X_t;\theta)\\
\text{(part. obsvd.)} \qquad U_t & = \psi_t(\pi_t;\theta)\
\end{align*}
where $\pi_t= {\sf P}(X_t\mid{\cal Y}_t)$ is the belief state; the use of the
common notation $\psi_t(\cdot;\theta)$ should not cause any confusion because the
two cases are treated in separate subsections.

The following operators are of interest (see \Sec{sec:fnspaces} for
definition of function space ${\cal M}$ and ${\cal U}$):
\begin{enumerate}
\item The operator $\Psi:{\cal M}\to {\cal U}$, defined as
\[
\Psi(\b) = \argmin_{U\in {\cal U}} {\sf J} (U;\b)
\]
\item The operator $\Xi: {\cal U} \to {\cal M}$, defined according to~\eqref{eq:rho} where $\beta$ is obtained using~\eqref{eq:mf-vars}.   
\end{enumerate}

Assuming that the two operators are well-defined, we have:

\medskip
\begin{definition}\label{def:MFE}
A \emph{mean-field equilibrium} (MFE) is a fixed point $\beta$ such that $\Xi(\Psi(\b)) = (\b)$.
\end{definition}


\subsection{Fully Observed Case}

For the fully observed problem,  a complete characterization of the
MFG equilibrium is described in the following
proposition.  The result is a minor extension of a similar result
appearing in our prior paper~\cite{olmez2021pre} (for
$\lambda^{\text{\tiny AR}}$ and $\eta=0$) and therefore its proof is
omitted.  

\begin{proposition}
The unique mean-field equilibrium is
\[
\b_t  = 0\quad \forall\;\;t\geq 0
\]
\end{proposition}

\medskip

This result represents a mathematical expression of the thought
experiment described in \Sec{sec:intro}.  Its main utility is to set up the problem whereby the effects of some of the
underlying assumptions -- perfect rationality and perfect information
-- can be investigated.  As with the preceding sections, the main
focus of the MFG modeling is in the partially observed settings, which is the
subject of the following subsection. 
 
\subsection{Partially Observed Case}
\label{sec:mfg_po}

Consider the space of probability
distributions on the belief space ${\cal P}^1\cup  {\cal P}^2$.  The
random vector $(S_t,A_t)$ is well-defined on the set $[t<\tau]$, and we
denote by $p_t(s,a)$ its density for $0\leq s \leq 1$, $s+a\leq1$:
\begin{align*}
&{\sf P}([s< S_t <  s+\ud s]\cap[a< A_t <  a+\ud a]\cap [t<\tau])\\
& = p_t(s,a)\ud a \ud s, \quad t\geq 0 
\end{align*}

\smallskip
\begin{assumption}\label{ass:ass3}
The density $p_t(s,1-s)=p_t(s,0)=0$ for all $t\geq 0$, $s \in [0,1]$. 
\end{assumption}
\smallskip


Under Assumption~\ref{ass:ass3} that an agent uses the optimal control
$U_t=U_t^{\text{opt}} = \psi_t(S_t,A_t)$ for $0\leq t<\tau$, the density process $\{p_t(s,a)\in [0,\infty):0\leq
s\leq 1,s+a \leq 1,t\geq 0\}$ solves the FPK equation whose derivation, along with the expressions for $f_\Sem (s,a,t)$ and $f_\Aem (s,a,t)$, appears in
Appendix~\ref{appdx:derivations}: 
\begin{align}
\frac{\partial p_t}{\partial t}(s,a)&=-\frac{\partial}{\partial s}\big(f_\Sem (s,a,t) p_t(s,a)\big) \nonumber \\
&-\frac{\partial}{\partial a}\big( f_\Aem (s,a,t)  p_t(s,a)\big)
-a\lambda^{\text{\tiny AI}}p_t(s,a)
\label{eq:FPK_POMFG_1}
\end{align}
where $p_0(s,a)$ is the initial density (assumed given).
By using the tower property,
\begin{align*}
\rho_t(\Aem) & = {\sf P}([X_t = \Aem]) \\
& = {\sf E}(\pi_t(\Aem)) = {\sf
  E}(A_t1_{[t<\tau]}) = \int_0^1 \int_0^{1-s} a p_t(s,a) \ud a \ud s
\end{align*}
and therefore we have
\begin{subequations}\label{eq: FPK_POMFG_2}
\begin{align}
\frac{\ud \rho_t}{\ud t} (\Iem)&= \lambda^{\text{\tiny AI}} \int_0^1 \int_0^{1-s} a p_t(s,a) \ud a \ud s -(\lambda^{\text{\tiny
                                 ID}} + \lambda^{\text{\tiny IR}}) \rho_t(\Iem)\\
\frac{\ud \rho_t}{\ud t}(x) &= ({\cal A}^\dagger \rho_t)(x),\quad x\in\{\Rem,\Dem\}
\end{align}
\end{subequations}
where expression for ${\cal A}^\dagger$ is obtained from the
transition graph.  

With a heterogeneous population, the notation $p_t(s,a;\theta)$ is used
to denote the
density conditioned on the attribute $\theta$.  The mean-field process is then consistently
obtained as 
\begin{align}
\b_t  = \sum_{\theta} {\sf p}(\theta) & \int_0^1 \int_0^{1-s} a \psi_t(s,a;\theta) p_t(s,a;\theta) \ud a \ud s
\label{eq:consistency_POMFG}
\end{align}

This completes the derivation of the system of equations for the
partially observed MFG: Eq.~\eqref{eq:FPK_POMFG_1}-\eqref{eq: FPK_POMFG_2} is the forward FPK equation.
Eq.~\eqref{eq:HJB_POMFG} is the backward HJB equation.
Eq.~\eqref{eq:consistency_POMFG} defines the consistency relationship
that links the two equations.  Its solution is an MFE (satisfies
Defn.~\ref{def:MFE}).  

The thrust of the ongoing work is to obtain numerical solutions of the
MFG equations and use it to explain differences between the data from
the omicron and delta waves.

\bibliographystyle{IEEEtran}
\bibliography{references,referencesnew}

\begin{thebibliography}{10}
\providecommand{\url}[1]{#1}
\csname url@samestyle\endcsname
\providecommand{\newblock}{\relax}
\providecommand{\bibinfo}[2]{#2}
\providecommand{\BIBentrySTDinterwordspacing}{\spaceskip=0pt\relax}
\providecommand{\BIBentryALTinterwordstretchfactor}{4}
\providecommand{\BIBentryALTinterwordspacing}{\spaceskip=\fontdimen2\font plus
\BIBentryALTinterwordstretchfactor\fontdimen3\font minus
  \fontdimen4\font\relax}
\providecommand{\BIBforeignlanguage}[2]{{%
\expandafter\ifx\csname l@#1\endcsname\relax
\typeout{** WARNING: IEEEtran.bst: No hyphenation pattern has been}%
\typeout{** loaded for the language `#1'. Using the pattern for}%
\typeout{** the default language instead.}%
\else
\language=\csname l@#1\endcsname
\fi
#2}}
\providecommand{\BIBdecl}{\relax}
\BIBdecl

\bibitem{olmez2021pre}
S.~Y. Olmez, S.~Aggarwal, J.~W. Kim, E.~Miehling, T.~Ba{\c{s}}ar, M.~West, and
  P.~G. Mehta, ``Modeling presymptomatic spread in epidemics via mean-field
  games,'' \emph{arXiv preprint arXiv:2111.10422}, 2021.

\bibitem{rivett2020screening}
L.~Rivett, S.~Sridhar, D.~Sparkes, M.~Routledge, N.~K. Jones, S.~Forrest,
  J.~Young, J.~Pereira-Dias, W.~L. Hamilton, M.~Ferris \emph{et~al.},
  ``Screening of healthcare workers for {SARS-CoV-2} highlights the role of
  asymptomatic carriage in {COVID-19} transmission,'' \emph{Elife}, vol.~9, p.
  e58728, 2020.

\bibitem{buitrago2020role}
D.~C. Buitrago-Garcia, D.~Egli-Gany, M.~J. Counotte, S.~Hossmann, H.~Imeri,
  G.~Salanti, and N.~Low, ``The role of asymptomatic {SARS-CoV-2} infections:
  rapid living systematic review and meta-analysis,'' \emph{MedRxiv}, 2020.

\bibitem{bender2021analysis}
J.~K. Bender, M.~Brandl, M.~H{\"o}hle, U.~Buchholz, and N.~Zeitlmann,
  ``Analysis of asymptomatic and presymptomatic transmission in {SARS-CoV-2}
  outbreak, {G}ermany, 2020,'' \emph{Emerging infectious diseases}, vol.~27,
  no.~4, p. 1159, 2021.

\bibitem{wei2020presymptomatic}
W.~E. Wei, Z.~Li, C.~J. Chiew, S.~E. Yong, M.~P. Toh, and V.~J. Lee,
  ``Presymptomatic transmission of {SARS-CoV-2—Singapore}, january 23--march
  16, 2020,'' \emph{Morbidity and Mortality Weekly Report}, vol.~69, no.~14, p.
  411, 2020.

\bibitem{elie2020contact}
R.~Elie, E.~Hubert, and G.~Turinici, ``Contact rate epidemic control of
  {COVID-19}: an equilibrium view,'' \emph{Mathematical Modelling of Natural
  Phenomena}, vol.~15, p.~35, 2020.

\bibitem{aurell2020optimal}
A.~Aurell, R.~Carmona, G.~Dayanikli, and M.~Lauriere, ``Optimal incentives to
  mitigate epidemics: a {S}tackelberg mean field game approach,'' \emph{arXiv
  preprint arXiv:2011.03105}, 2020.

\bibitem{aurell2021finite}
------, ``Finite state graphon games with applications to epidemics,''
  \emph{arXiv preprint arXiv:2106.07859}, 2021.

\bibitem{hubert2020incentives}
E.~Hubert, T.~Mastrolia, D.~Possama{\"\i}, and X.~Warin, ``Incentives,
  lockdown, and testing: from thucydides's analysis to the covid-19 pandemic,''
  \emph{arXiv preprint arXiv:2009.00484}, 2020.

\bibitem{lee2020controlling}
W.~Lee, S.~Liu, H.~Tembine, W.~Li, and S.~Osher, ``Controlling propagation of
  epidemics via mean-field control,'' \emph{arXiv preprint arXiv:2006.01249},
  2020.

\bibitem{cho2020mean}
S.~Cho, ``Mean-field game analysis of {SIR} model with social distancing,''
  \emph{arXiv preprint arXiv:2005.06758}, 2020.

\bibitem{doncel2020mean}
J.~Doncel, N.~Gast, and B.~Gaujal, ``A mean field game analysis of {SIR}
  dynamics with vaccination,'' \emph{Probability in the Engineering and
  Informational Sciences}, pp. 1--18, 2020.

\bibitem{tembine2020covid}
H.~Tembine, ``{COVID-19}: Data-driven mean-field-type game perspective,''
  \emph{Games}, vol.~11, no.~4, p.~51, 2020.

\bibitem{saldi2019partially}
N.~Saldi, T.~Ba{\c{s}}ar, and M.~Raginsky, ``Partially-observed discrete-time
  risk-sensitive mean-field games,'' in \emph{2019 IEEE 58th Conference on
  Decision and Control (CDC)}.\hskip 1em plus 0.5em minus 0.4em\relax IEEE,
  2019, pp. 317--322.

\bibitem{saldi2019po}
\BIBentryALTinterwordspacing
N.~Saldi, T.~Başar, and M.~Raginsky, ``Approximate nash equilibria in
  partially observed stochastic games with mean-field interactions,''
  \emph{Mathematics of Operations Research}, vol.~44, no.~3, pp. 1006--1033,
  2019. [Online]. Available: \url{https://doi.org/10.1287/moor.2018.0957}
\BIBentrySTDinterwordspacing

\bibitem{confortola2013filtering}
F.~Confortola and M.~Fuhrman, ``Filtering of continuous-time markov chains with
  noise-free observation and applications,'' \emph{Stochastics An International
  Journal of Probability and Stochastic Processes}, vol.~85, no.~2, pp.
  216--251, 2013.

\bibitem{nyt2022where}
\BIBentryALTinterwordspacing
S.~Cobey, J.~Bloom, T.~Starr, and N.~Lash, ``We study virus evolution. here's
  where we think the coronavirus is going.'' \emph{New York Times}, 2022.
  [Online]. Available:
  \url{https://www.nytimes.com/interactive/2022/03/28/opinion/coronavirus-mutation-future.html}
\BIBentrySTDinterwordspacing

\bibitem{cassidy09rational}
J.~Cassidy, ``Rational irrationality: The real reason that capitalism is so
  crash-prone,'' Oct. 5th 2009.

\bibitem{Yin_synch_game}
H.~Yin, P.~G. Mehta, S.~P. Meyn, and U.~V. Shanbhag, ``Synchronization of
  oscillators is a game,'' in \emph{IEEE Transactions on Automatic Control},
  vol.~57, no.~4, April 2012, pp. 920--935.

\bibitem{xiong2008book}
J.~Xiong, \emph{An Introduction to Stochastic Filtering Theory}.\hskip 1em plus
  0.5em minus 0.4em\relax Oxford University Press, 2008.

\bibitem{fleming2006controlled}
W.~H. Fleming and H.~M. Soner, \emph{\BIBforeignlanguage{eng}{Controlled Markov
  Processes and Viscosity Solutions}}, 2nd~ed., ser. Stochastic Modelling and
  Applied Probability ; 25.\hskip 1em plus 0.5em minus 0.4em\relax New York:
  Springer-Verlag, 2006.

\end{thebibliography}

\appendix

\subsection{Proof of Proposition~\ref{prop:sta_sol}}
With $\b_t = \bar{\b}$, the HJB equation for $v_t(\Sem)$ is 
\begin{align*}
-\frac{\ud v_t}{\ud t} (\Sem)& + (\gamma+\eta) v_t (\Sem) \\
&= \eta \bar{\phi} (\Aem) + \min_{u \in [0,1]} \left( \beb \left( \bar{\phi} (\Aem) - v_t (\Sem) \right) - \v \right) u
\end{align*}
We investigate its stationary solutions $v_t(\Sem) =
\bar{\phi}(\Sem)$ in which case the stationary HJB equation is
\begin{align}
(\gamma+\eta) \bar{\phi} (\Sem) &= \eta \bar{\phi} (\Aem) \nonumber \\
&+ \min_{u \in [0,1]} \underbrace{ \left( \beb \left( \bar{\phi} (\Aem) - \bar{\phi} (\Sem) \right) - \v \right) }_{=:M} u \label{eq:stat_HJB_fully_obs_phi}
\end{align}
We have the  following two cases:
\begin{itemize}
\item If $\bar{\beta}<\frac{\v}{\lambda^{\text{\tiny
                                SA}}\f(\Aem)} \left( 1+\frac{\eta}{\gamma} \right)$, then $\bar{\phi} (\Sem)= \frac{(\lambda^{\text{\tiny SA}} \bar{\beta}+\eta)\bar{\phi} (\Aem)- \v}{\lambda^{\text{\tiny SA}}\bar{\b} + \gamma + \eta}$ solves the HJB equation with the minimizing choice of
$u=1$, because
\begin{align*}
&M
        = \beb \left( \frac{\gamma \bar{\phi} (\Aem) + \v}{\beb + \gamma + \eta} \right) - \v<0\\
\iff& \beb \v + \gamma \v + \eta \v > \beb \gamma \bar{\phi} (\Aem) + \beb \v\\
\iff & \bar{\beta}<\frac{\v}{\lambda^{\text{\tiny
                                SA}}\f(\Aem)} \left( 1+\frac{\eta}{\gamma} \right)
\end{align*}
\item If $\bar{\beta} > \frac{\v}{\lambda^{\text{\tiny
                                SA}}\f(\Aem)} \left( 1+\frac{\eta}{\gamma} \right)$, then $\bar{\phi}(\Sem) =\frac{\eta}{\gamma+\eta} \bar{\phi} (\Aem)$  solves the HJB equation with the minimizing choice of
$u=0$, because
\[
M = \frac{\beb \gamma \bar{\phi} (\Aem)}{\gamma+\eta} - \v>0
\]
\end{itemize}

\subsection{Proof of Proposition~\ref{prop:filter}}\label{appdx:filter}

Because $A_0=0$, and the solution is continuous as a function of time
$t$, it suffices to show that 
\begin{align*}
A_t = \frac{\beb + \eta}{\la} \implies \frac{\ud A_t}{\ud t} < 0
\end{align*}
Setting $A_t = \frac{\beb + \eta}{\la}$ on the righthand side
of~\eqref{eq:filter}, 
\begin{align*}
\frac{\ud A_t}{\ud t} = (\beb + \eta) \left( S_t -1 + \frac{\beb+\eta-\lar}{\la} \right)
\end{align*}
The desired conclusion follows because $S_t+ A_t \leq 1$ and therefore
$S_t - 1 +\frac{\beb +
  \eta}{\la} \leq 0$ if $A_t = \frac{\beb + \eta}{\la}$. 

\subsection{Derivations}\label{appdx:derivations}
Let ${\cal P}:={\cal P}^1\cup {\cal P}^2$.  
The belief process $\pi:=\{\pi_t\in{\cal P}:t\geq
0\}$ is a Markov process~\cite[Theorem 1.7]{xiong2008book}.  The HJB equations and FPK equations
are easily derived once we obtain the infinitesimal generator of the
process.  

\smallskip

\newP{Infinitesimal generator} Consider a smooth test function $v:{\cal P} \to \Re$. 
Let $\mu\in {\cal P} $.  The infinitesimal generator
(for the general time-inhomogeneous case) is
	\[
	({\cal A}_t^u v)(\mu) = \lim_{\delta t\downarrow 0}
        \frac{\E(v(\pi_{t+\delta t})|\pi_t=\mu)-v(\mu)}{\delta t}
	\]
There are two cases to consider:

\smallskip

\noindent $\bullet$ If $\mu\in {\cal P}^2$, upon identifying the
measures $\{\delta_\Iem,\delta_\Rem,\delta_\Dem\}$ with the states
$\{\Iem,\Rem,\Dem\}$, the generator is the same as the generator for
the Markov process.  

\smallskip

\noindent $\bullet$
If $\mu\in{\cal P}^1 \setminus \{ \delta_\Rem \}$, then using the coordinates $(s,a)$ for ${\cal
    P}^1$, $\mu = [s,a,0,1-s-a,0]$ for $a\in[0,1]$.  With $\pi_t=\mu$,
  in the asymptotic limit as $\delta t\to 0$,
\[
\pi_{t+\delta t} = \begin{cases}
	\mu + e_{\Sem} f_{\Sem}(s,a,t) \delta t \\+ e_{\Aem} f_{\Aem}(s,a,t) \delta t + o(\delta t) & \text{w.p.} \;\;
        (1-a\lambda^{\text{\tiny AI}}\delta t) + o(\delta t)\\
	\delta_\Iem &\text{w.p.} \;\;
        a\lambda^{\text{\tiny AI}}\delta t + o(\delta t)
	\end{cases}
	\]   
where 
\begin{align*}
f_{\Sem}(s,a,t)&=\left( -\be u - \eta + a \la \right) s \\
f_\Aem (s,a,t) &= \left( \be u + \eta \right) s + a \left( -\la - \lar + \la a \right)
\end{align*}
and $e_\Sem = [1,0,0,-1,0]$, $e_\Aem = [0,1,0,-1,0]$.
Denoting $v(\mu) = \phi(s,a)$, the generator is then easily calculated to be
\begin{align*}
({\cal A}_t^u v)(\mu) &= f_\Sem(s,a,t)  \frac{\partial \phi}{\partial s}(s,a)+f_\Aem(s,a,t)  \frac{\partial \phi}{\partial a}(s,a)\\& +
a \lambda^{\text{\tiny AI}} (v(\delta_\Iem)-\phi(s,a) )
\end{align*}
where the superscript $u$ denotes the fact that $f_\Sem(s,a,t)$, $f_\Aem(s,a,t)$, and therefore 
also the generator, depends also upon
$u$.  The subscript $t$ denotes the fact that the generator is for a
time-inhomogeneous Markov process (because $\beta$ may depend
upon time).  

\smallskip

\newP{Derivation of the HJB equation} 
The HJB equation (See Sec.III.7 of~\cite{fleming2006controlled}) is
\[
- \frac{\partial v_t}{\partial t} (\mu) + \gamma v_t(\mu) = \min_{u\in[0,1]} \left( {\cal A}_t^u v_t(\mu) + \mu
  (c(\cdot,u;\v)) \right)
\]
In order to calculate $\mu\big(c(\cdot,u;\v)\big)$, we need to assume a running cost for state $\Rem$. If we assume $c(\Rem,\cdot;\cdot) = \gamma \bar{\phi} (\Rem)$, then we get
\begin{align*}
\int_0^\infty c(\Rem,\cdot;\cdot) e^{-\gamma t} dt = \bar{\phi} (\Rem)
\end{align*}
as desired. Let
\begin{align*}
&\mathcal{L} \phi := (1-s-a) \gamma \bar{\phi} (\Rem) + a \la \left( \bar{\phi} (\Iem) - \phi \right)\\
&+(a\la - \eta) s \frac{\partial \phi}{\partial s} + \left( \eta s + a (a \la - \la -\lar) \right) \frac{\partial \phi}{\partial a}\\
&{\cal M}_t(s,a) := a-(s+a) \alpha + \be s \left( \frac{\partial \phi_t}{\partial a} - \frac{\partial \phi_t}{\partial s} \right)
\end{align*}
Thus for $\mu = [s,a,0,1-s-a,0] $ and $v_t(\mu) = \phi_t(s,a)$, the HJB
equation~\eqref{eq:HJB_POMFG} is obtained because
\begin{align*}
\mu\big(c(\cdot,u;\v)\big) &= (1-s-a)\gamma \bar{\phi} (\Rem) + s(-\v u) + a (1-\v) u\\
&= (1-s-a) \gamma \bar{\phi} (\Rem) + (a-(s+a)\v) u
\end{align*}
The boundary conditions~\eqref{eq:bcs} follow from the observation that
\begin{align*}
&\phi_t(s,a) \\
&= s \min_{U \in L_{\mathcal{F}}^2} {\sf E} \bigg( \int_t^T e^{-\gamma s} c(X_s,U_s;\v) ds + e^{-\gamma T} \phi(X_T)\\ &\mid X_t = \Sem \bigg)
 + a \bar{\phi} (\Aem) + (1-s-a) \bar{\phi} (\Rem)
\end{align*}

\smallskip

\newP{Derivation of the FPK equation}
We derive the adjoint of the generator ${\cal A}^u$ where dependence
on $t$ is suppressed for notational ease. Let $\rho$ be a measure on ${\cal P}$. On ${\cal P}^1$, $\rho$ has density $p(s,a)$. Consider $\rho({\cal A}^uv) = \int {\cal A}^uv(\mu) \rho(\ud \mu)=$ 
\begin{align}
&\int_0^1 \int_0^{1-s} p(s,a) \bigg[ f_\Sem(s,a) \frac{\partial \phi}{\partial s} + f_\Aem(s,a) \frac{\partial \phi}{\partial a} \nonumber \\
&+a\la \big( v(\delta_\Iem) - \phi(s,a) \big)   \bigg] \ud a \ud s - \rho(\delta_\Iem) (\lr+\ld) v(\delta_\Iem) \nonumber \\
& +\rho(\delta_\Rem) \lr v(\delta_\Rem) + \rho(\delta_\Dem) \ld v(\delta_\Dem)
\label{eq:fpk_der_1}
\end{align}
First note that
\begin{align*}
&\int_0^{1-a} p(s,a) f_\Sem (s,a) \frac{\partial \phi}{\partial s} \ud s = p(s,a) f_\Sem (s,a) \phi (s,a) \vert_{s=0}^{s=1-a}\\
&-\int_0^{1-a} \phi(s,a) \frac{\partial}{\partial s} \big( p(s,a) f_\Sem (s,a) \big) \ud s 
\end{align*}
where boundary terms vanish because $f_\Sem(0,a)=0$ and $p(1-a,a)=0$
(Assumption~\ref{ass:ass3}). One can also write
\begin{align*}
&\int_0^{1-s} p(s,a) f_\Aem (s,a) \frac{\partial \phi}{\partial a} \ud a = p(s,a) f_\Aem (s,a) \phi (s,a) \vert_{a=0}^{a=1-s}\\
&-\int_0^{1-s} \phi(s,a) \frac{\partial}{\partial a} \big( p(s,a) f_\Aem (s,a) \big) \ud a 
\end{align*}
where boundary terms vanish because $p(s,1-s)=p(s,0)=0$ (Assumption~\ref{ass:ass3}). Using these two results, one can rewrite the righthand side of~\eqref{eq:fpk_der_1} as
\begin{align*}
&-\int_0^1 \int_0^{1-s} \phi(s,a) \bigg[ \frac{\partial}{\partial s} \big( p(s,a) f_\Sem (s,a) \big) + \frac{\partial}{\partial a} \big( p(s,a) f_\Aem (s,a) \big) \\
& +a\la p(s,a) \bigg] \ud a \ud s +v(\delta_\Rem) \lr \rho(\delta_\Rem) + v(\delta_\Dem) \ld \rho(\delta_\Dem)\\
& + v(\delta_\Iem) \left( \int_0^1 \int_0^{1-s} a \la p(s,a) \ud a \ud s - (\lr+\ld) \rho(\delta_\Iem) \right)
\end{align*}

\subsection{Numerical Algorithm}
\label{sec:numol}

The HJB equation~\eqref{eq:HJB_POMFG} is numerically solved using the method of
lines.  In this method, the spatial coordinate $(s,a)\in D$ is
discretized over a finite grid.  The partial derivatives $\partial
\phi_t / \partial s$ and $\partial \phi_t / \partial a$ are
approximated using a finite difference approximation at each point on
the grid.  The time variable $t$ is treated as a continuous variable
and after spatial discretization, the resulting ordinary differential
equations are numerically integrated (backward in time) using a
standard numerical procedure.  For this purpose, a sufficiently large
terminal time $T$ is chosen with the following terminal condition:
\begin{align*}
\phi_T(s,a) = \begin{cases}
a \bar{\phi} (\Aem) + (1-a) \bar{\phi} (\Rem) & a \leq \frac{1-s}{2} \\
\phi_T(s,1-s-a) \\
+ (2a+s-1) \left( \bar{\phi} (\Aem) - \bar{\phi} (\Rem) \right) & a > \frac{1-s}{2}
\end{cases}
\label{eq:arb_tc}
\end{align*}
The algorithm was used to compute the stationary solutions depicted in
Fig.~\ref{fig:u_sta_beta}.

\end{document}